\documentclass{amsart}

\usepackage{amsmath, amssymb, mathrsfs, dsfont}

\numberwithin{equation}{section}

\begin{document}

\title{Characterization of $\gamma$-factors: the Asai case}


\author{Guy Henniart}
\address{Guy henniart} %\\ Institut Universitaire de France et Univ. Paris-Sud \\ Laboratoire de Math\'ematiques d'Orsay \\ CNRS \\ Orsay cedex F-91405, France}
\email{Guy.Henniart@math.u-psud.fr}

\author{Luis Lomel\'i}
\address{Luis Lomel\'i} % \\ }
\email{lomeli@math.purdue.edu}
%\email{llomeli@math.ou.edu}


\date{January 2012}

\begin{abstract}
Let $E$ be a separable quadratic extension of a locally compact field $F$ of positive characteristic. Asai $\gamma$-factors are defined for smooth irreducible representations $\pi$ of ${\rm GL}_n(E)$. If $\sigma$ is the Weil-Deligne representation of $\mathcal{W}_E$ corresponding to $\pi$ under the local Langlands correspondence, we show that the Asai $\gamma$-factor is the same as the Deligne-Langlands $\gamma$-factor of the Weil-Deligne representation of $\mathcal{W}_F$ obtained from $\sigma$ under tensor induction. This is achieved by proving that Asai $\gamma$-factors are characterized by their local properties together with their role in global functional equations for $L$-functions. As an immediate application, we establish the stability property of $\gamma$-factors under twists by highly ramified characters.
\end{abstract}

\maketitle

\section{Introduction}

Let $F$ be a locally compact field of characteristic $p$. Let $\psi$ be a non-trivial character of $F$ and $\pi$ a smooth irreducible representation of ${\rm GL}_n(F)$, where $n$ is a positive integer. Let $r$ denote either ${\rm Sym}^2$ or $\wedge^2$. In particular, if $\rho_n$ denotes the standard representation of ${\rm GL}_n(\mathbb{C})$, then $r$ denotes either ${\rm Sym}^2 \rho_n$ or $\wedge^2 \rho_n$. If $\tau$ is a Weil-Deligne representation, let $\gamma_F^{\rm Gal}(s,\tau,\psi)$, $s \in \mathbb{C}$, be the $\gamma$-factor defined by Deligne and Langlands \cite{d,t}. In \cite{hl}, the authors establish the equality of $\gamma$-factors:
\begin{equation*}
   \gamma(s,\pi,r,\psi) = \gamma_F^{\rm Gal}(s,r \circ \sigma,\psi),
\end{equation*}
where the factor on the left is defined in \cite{luis} and $\sigma$ on the right is the transfer of $\pi$ to a Weil-Deligne representation under the local Langlands correspondence \cite{lrs}.

In this paper, we address the case of Asai $\gamma$-factors and related $L$- and $\varepsilon$-factors. Asai factors can be seen as a generalization of those studied in \cite{asai} by T. Asai. Let $E/F$ be a separable quadratic extension with Galois group ${\rm Gal}(E/F) = \left\{ 1, \theta \right\}$. Fix a separable algebraic closure $\overline{F}$ of $F$ containing $E$. Let $\pi$ be a smooth irreducible representation of ${\rm GL}_n(E)$. The $L$-group of ${\rm Res}_{E/F}{\rm GL}_n$ is ${\rm GL}_n(\mathbb{C}) \times {\rm GL}_n(\mathbb{C}) \rtimes \mathcal{W}_F$, where the Weil group $\mathcal{W}_F$ of $\overline{F}/F$ acts via the Galois group and $\theta$ permutes the two copies of ${\rm GL}_n(\mathbb{C})$. Its Asai representation $r_{\mathcal{A}}$ can be defined for every $n$ by
\begin{align*}
   r_{\mathcal{A}} & : {\rm GL}_n(\mathbb{C}) \times {\rm GL}_n(\mathbb{C}) \rtimes {\rm Gal}(E/F) \rightarrow {\rm GL}_{n^2}(\mathbb{C}), \\
   r_{\mathcal{A}} & (x,y,1) = (x \otimes y) \text{ and } r_{\mathcal{A}}(x,y,\theta) = (y \otimes x).
\end{align*}
Asai $\gamma$-factors $\gamma_{E/F}(s,\pi,r_{\mathcal{A}},\psi)$ are defined in \cite{luis} for generic $\pi$ by developing the Langlands-Shahidi method in positive characteristic. Asai local factors are defined for a general smooth irreducible $\pi$ using Langlands' classification in \S~4.2. Let $\sigma$ be the transfer of $\pi$ to a Weil-Deligne representation of $\mathcal{W}_E$ under the local Langlands correspondence. We prove the equality
\begin{equation} \label{maineq}
   \gamma_{E/F}(s,\pi,r_{\mathcal{A}},\psi) = \gamma_F^{\rm Gal}(s,{}^\otimes{\rm I}(\sigma),\psi),
\end{equation}
where ${}^\otimes{\rm I}(\sigma)$ denotes the representation of $\mathcal{W}_F$ obtained from $\sigma$ by tensor induction. In the case of characteristic zero, equation~(1.1) for $n=2$ is known \cite{krish,r} and up to roots of unity in general \cite{he10}.

Theorem~3.3 is proved via a characterization of Asai $\gamma$-factors involving local properties together with their connection with the global theory by means of a functional equation described below~(1.2). More precisely, the local properties of $\gamma_{E/F}(s,\pi,r_{\mathcal{A}},\psi)$ include: a naturality property with respect to isomorphisms of quadratic extensions $E/F$; an isomorphism property pertaining to $\pi$; a dependence on the additive character $\psi$, which can be made explicit; a crucial multiplicativity property with respect to parabolic induction, which reflects the influence of taking tensor induction on a direct sum of Weil-Deligne representations; and finally, (1.1), is needed whenever the representation $\pi$ is the generic component of an unramified principal series (this can be reduced, by multiplicativity, to the case where $n=1$).

As indicated above, part of our argument is global. Now, let $K/k$ be a separable quadratic extension of global function fields of characteristic $p$. Let $\Psi$ be a non-trivial character of $k \backslash \mathbb{A}_k$ and let $\Pi = \otimes' \Pi_v$ be a cuspidal automorphic representation of ${\rm GL}_n(\mathbb{A}_K)$. If a place $v$ of $k$ is split, then $\gamma_{K_v/k_v}(s,\Pi_v,r_{\mathcal{A}},\Psi_v) = \gamma_{k_v}(s,\pi_1 \times \pi_2,\psi)$, a Rankin-Selberg $\gamma$-factor (see \S~6.2 of \cite{luis} and the appendix). Let $S$ be a finite set of places such that $K/k$, $\Pi$ and $\Psi$ are unramified for $v \notin S$. The global functional equation is given by
\begin{equation}
   L^S(s,\Pi,r_{\mathcal{A}}) = \prod_{v \in S} \gamma_{K_v/k_v}(s,\Pi_v,r_{\mathcal{A}},\Psi_v) L^S(1-s,\widetilde{\Pi},r_{\mathcal{A}}),
\end{equation}
where
\begin{equation*}
   L^S(s,\Pi,r_{\mathcal{A}}) = \prod_{v \notin S} L(s,\Pi_v,r_{\mathcal{A}}).
\end{equation*}

The approach to prove the main result is that of \cite{hl} where equality of local $\gamma$-factors is obtained from the global functional equation for special cases of cuspidal automorphic representations where there is control over the places of $k$ where ramification can occur. In fact, we directly establish a local-to-global argument for the case of a cuspidal, tamely ramified representation $\pi$ of ${\rm GL}_n(E)$ (hence $\pi$ of level zero) via the Grundwald-Wang theorem. Then, the general problem is reduced to the case of a tamely ramified representation $\pi$. This is done by using a local-to-global result due to Gabber and Katz for $\ell$-adic representations of the Galois group \cite{katz}, and translating it via the global Langlands correspondence \cite{laff}.

In the Langlands-Shahidi method, $\pi$ is assumed to be generic. But, using the Langlands-Zelevinsky classification together with multiplicativity it is possible to work with a general irreducible smooth representation. See \S~4.2, where the definition of local Asai $L$-functions and root numbers, which are obtained via Asai $\gamma$-factors, is extended to the general case. In Theorem~4.3, we show in this general setting that the local $L$- and $\varepsilon$-factors are the same as the corresponding Galois factors. In \S~4.5 we take the opportunity to write down a stability property of $\gamma$-factors that is not known in characteristic zero.

The authors would like to thank F. Shahidi and B. Gross for mathematical communications. The second author also thanks O. Gabber, M. Krishnamurthy, P. Kutzko, L. Lafforgue, F. Shahidi and Y. Ye for useful mathematical discussions. He also thanks the University of Iowa, Purdue University and the IH\'ES  for their hospitality and support while the article was written.

\section{Asai $\gamma$-factors and tensor induction}

\noindent{\bf 2.1} Fix a prime $p$ and define $\mathscr{A}_{\rm quad}(p)$ to be the class whose objects are triples $(E/F,\pi,\psi)$ consisting of: a separable quadratic extension $E/F$ of locally compact fields of characteristic $p$; a smooth irreducible representation $\pi$ of ${\rm GL}_n(E)$, for some positive integer $n$; and, a non-trivial character $\psi$ of $F$.

Given a triple $(E/F,\pi,\psi) \in \mathscr{A}_{\rm quad}(p)$, we say it is of degree $n$ if $\pi$ is a representation of ${\rm GL}_n(E)$. Also, we say $(E/F,\pi,\psi) \in \mathscr{A}_{\rm quad}(p)$ is generic if $\pi$ is generic. Let $x \mapsto \bar{x}$ denote conjugation in $E/F$, i.e., the non-trivial automorphism of $E/F$. Let $q$ be the cardinality of the residue field of $F$. Let $\pi^{\rm conj}$ denote the representation obtained from $\pi$ by conjugation in $E/F$ and let $\widetilde{\pi}$ denote the contragredient of $\pi$.

Asai $\gamma$-factors give a rule which associates a rational function $\gamma_{E/F}(s,\pi,r_{\mathcal{A}},\psi) \in \mathbb{C}(q^{-s})$ to each triple $(E/F,\pi,\psi) \in \mathscr{A}_{\rm quad}(p)$. They are initially defined in \cite{luis} when $(E/F,\pi,\psi) \in \mathscr{A}_{\rm quad}(p)$ is generic. This definition is extended to all of $\mathscr{A}_{\rm quad}(p)$ in \S~4.2 by means of Langlands classification for $\mathfrak{p}$-adic reductive groups.

\bigskip

\noindent{\bf 2.2} Given a locally compact non-archimedean field $F$ of characteristic $p$ and a separable quadratic extension $E/F$, fix a separable algebraic closure $\bar{F}$ containing $E$. Let $\mathscr{G}_{\rm quad}(p)$ be the class whose objects are triples $(E/F,\sigma,\psi)$ consisting of: a separable quadratic extension $E/F$ of locally compact fields; an $n$-dimensional Frobenius semisimple Weil-Deligne representation of the Weil group $\mathcal{W}_E$ of $E$, for some positive integer $n$; and, a non-trivial character $\psi$ of $F$.

The $L$-functions and $\varepsilon$-factors attached to a Weil-Deligne representation $\tau$ of $\mathcal{W}_F$, give rise to the Galois $\gamma$-factor
\begin{equation*}
   \gamma_F^{\rm Gal}(s,\tau,\psi) = \varepsilon(s,\tau,\psi) L(1-s,\widetilde{\tau})/L(s,\tau).
\end{equation*}
Here, $\widetilde{\tau}$ denotes the contragredient of $\tau$. This Galois $\gamma$-factor depends only on the isomorphism class of $\tau$ and $\psi$.

Given a representation $\sigma$ of $\mathcal{W}_E$, let ${}^\otimes{\rm I}(\sigma)$ be the representation of $\mathcal{W}_F$ obtained from $\sigma$ via tensor induction (see \S~13 of \cite{cr}). For each triple $(E/F,\sigma,\psi) \in \mathscr{G}_{\rm quad.}(p)$, the Galois $\gamma$-factor
\begin{equation*}
   \gamma_F^{\rm Gal}(s,{}^\otimes{\rm I}(\sigma),\psi)
\end{equation*}
is a rational function on $q^{-s}$ that arises in connection with Asai $\gamma$-factors. These factors satisfy a number of known properties, including a multiplicativity property reflecting the decomposition rule:
\begin{equation*}
   {}^\otimes{\rm I}(\sigma_1 \oplus \sigma_2) \simeq {}^\otimes{\rm I}(\sigma_1) \oplus {}^\otimes{\rm I}(\sigma_2) \oplus {\rm Ind}_E^F(\sigma_1 \otimes \sigma_2^{\rm conj}),
\end{equation*}
where $\sigma_2^{\rm conj}$ is obtained from $\sigma_2$ via conjugation in $E/F$. See \S~3 for the corresponding properties of Asai $\gamma$-factors.

\bigskip

\noindent{\bf 2.3. Remark.} The reader should consult \S~7 of \cite{ggp} for more details on $L$-groups of classical groups and their representations. Let $\eta_{E/F}$ be the character of $F^\times$ defining $E$ via local class field theory. Our representation $r_{\mathcal{A}}$ is ${\rm As}^+$, in the notation of \cite{ggp}, while ${\rm As}^-$ can be obtained from $r_{\mathcal{A}}$ after twisting by $\eta_{E/F}$. The latter representation arises in the Langlands-Shahidi method when considering odd unitary groups (see \S~6.4 of \cite{luis}). Also, twists of $r_{\mathcal{A}}$ by characters arise by considering the group of similitudes corresponding to the hermitian form $h(x,y)$ of \S~4.1.

\section{Characterization of Asai factors}

\noindent{\bf 3.1} We first list the local properties of Asai $\gamma$-factors which are used in the characterization. The existence of a system of $\gamma$-factors on $\mathscr{A}(p)$ is provided in \S~6 of \cite{luis} for generic representations. We remark on the dependence on the character $\psi$ in \S~4.1.

\begin{itemize}

\item[(i)] (Naturality). \emph{Let $(E/F,\pi,\psi) \in \mathscr{A}_{\rm quad}(p)$ be of degree $n$, and let $\eta$ be an isomorphism $\eta : E'/F' \simeq E/F$, more precisely, an isomorphism $\eta: E' \rightarrow E$ of local fields such that $\eta \vert_{F'}$ maps $F'$ into $F$. Then $\psi' = \psi \circ \eta \vert_{F'}$ is a non-trivial character of $F'$. Also, via $\eta$, $\pi$ defines a smooth irreducible representation $\pi'$ of ${\rm GL}_n(E')$. Then
\begin{equation*}
   \gamma_{E/F}(s,\pi,r_{\mathcal{A}},\psi) = \gamma_{E'/F'}(s,\pi',r_{\mathcal{A}},\psi').
\end{equation*}
}

\item[(ii)] (Isomorphism). \emph{Let $(E/F,\pi,\psi) \in \mathscr{A}_{\rm quad}(p)$ be of degree $n$, and let $\pi'$ be a smooth irreducible representation of ${\rm GL}_n(E)$ isomorphic to $\pi$. Then
\begin{equation*}
   \gamma_{E/F}(s,\pi',r_{\mathcal{A}},\psi) = \gamma_{E/F}(s,\pi,r_{\mathcal{A}},\psi).
\end{equation*}
}

\item[(iii)] (Compatibility with class field theory). \emph{A triple $(E/F,\chi,\psi)$ in $\mathscr{A}_{\rm quad}(p)$ of degree $n=1$, gives a character $\chi \vert_{F^\times}$ of ${\rm GL}_1(F)$. Then
\begin{equation*}
   \gamma_{E/F}(s,\chi,r_{\mathcal{A}},\psi) = \gamma_{F}(s,\chi \vert_{F^\times},\psi),
\end{equation*}
where the factor on the right hand side is an abelian $\gamma$-factor of Tate's thesis.
}

\item[(iv)] (Dependence on $\psi$). \emph{Let $(E/F,\pi,\psi) \in \mathscr{A}_{\rm quad}(p)$ be of degree $n$, and let $a \in F^\times$. Then $\psi^a : F \rightarrow \mathbb{C}^\times$, $x \mapsto \psi(ax)$, is a non-trivial character and we have
\begin{equation*}
   \gamma_{E/F}(s,\pi,r_{\mathcal{A}},\psi^a) = \omega_\pi(a)^n \left| a \right|_F^{n^2(s-\frac{1}{2})} \gamma_{E/F}(s,\pi,r_{\mathcal{A}},\psi).
\end{equation*}
}

\item[(v)] (Multiplicativity). \emph{Let $(E/F,\pi,\psi) \in \mathscr{A}_{\rm quad}(p)$. For $i=1, \ldots, d$, let $(E/F,\pi_i,\psi) \in \mathscr{A}_{\rm quad}(p)$ be of degree $n_i$ with each $\pi_i$ quasi-tempered and $n = n_1 + \cdots + n_d$. Assume that $\pi$ is the Langlands quotient of the representation of ${\rm GL}_n(E)$ obtained via unitary parabolic induction from $\pi_1 \otimes \cdots \otimes \pi_d$. Then
\begin{equation*}
   \gamma_{E/F}(s,\pi,r_{\mathcal{A}},\psi) = \prod_{i=1}^d \gamma_{E/F}(s,\pi_i,r_{\mathcal{A}},\psi) \prod_{i<j} \gamma_E(s,\pi_i \times \pi_j^{\rm conj},\psi \circ {\rm Tr}_{E/F}).
\end{equation*}
}

\item[(vi)] (Twists by unramified characters). \emph{Let $(E/F,\pi,\psi) \in \mathscr{A}_{\rm quad}(p)$. Then
\begin{equation*}
   \gamma_{\mathcal{A}}(s,\pi \otimes \left| \det(\cdot) \right|_E^{\frac{s_0}{2}},\psi) = \gamma_{\mathcal{A}}(s + s_0,\pi,\psi).
\end{equation*}
}

\end{itemize}

\noindent{\bf 3.2} The link to the global theory is provided by the following property (see \S~5 and 6.4(vi) of \cite{luis}):

\begin{itemize}

\item[(vii)] (Global functional equation). \emph{Let $K/k$ be a separable quadratic extension of global function fields of characteristic $p$. Let $\Psi$ be a non-trivial character of $\mathbb{A}_k/k$ and let $\Pi = \otimes' \Pi_v$ be an automorphic cuspidal representation of ${\rm GL}_n(\mathbb{A}_K)$. Given a place $v$ of $k$, let $K_v = K \otimes k_v$. If $K_v \simeq k_v \times k_v$, then $\Pi_v \simeq \pi_1 \otimes \pi_2$, where $\pi_1$ and $\pi_2$ are smooth irreducible representations of ${\rm GL}_n(k_v)$ and
\begin{equation*}
   \gamma_{K_v/k_v}(s,\Pi_v,r_{\mathcal{A}},\Psi_v) = \gamma_{k_v}(s,\pi_1 \times \pi_2,\psi).
\end{equation*}
Let $S$ be a finite set of places such that $K/k$, $\Pi$ and $\Psi$ are unramified outside of $S$. Then
\begin{equation*}
   L^{S}(s,\Pi,r_{\mathcal{A}}) = \prod_{v \in S} \gamma_{K_v/k_v}(s,\Pi_v,r_{\mathcal{A}},\Psi_v) L^S(1-s,\widetilde{\Pi},r_{\mathcal{A}}),
\end{equation*}
where
\begin{equation*}
   L^S(s,\Pi,r_{\mathcal{A}}) = \prod_{v \notin S} L(s,\Pi_v,r_{\mathcal{A}}).
\end{equation*}
}

\end{itemize}

\bigskip

\noindent{\bf 3.3. Theorem.} \emph{Let $\gamma_{\mathcal{A}}$ be a rule which assigns a rational function $\gamma_{\mathcal{A}}(s,\pi,\psi) \in \mathbb{C}(q^{-s})$ to each triple $(E/F,\pi,\psi) \in \mathscr{A}_{\rm quad}(p)$ satisfying properties (i)--(vi). If $(E/F,\pi,\psi) \in \mathscr{A}_{\rm quad}(p)$, let $(E/F,\sigma,\psi) \in \mathscr{G}_{\rm quad}(p)$ be the triple obtained by taking $\sigma = \sigma(\pi)$ to be the transfer of $\pi$ under the local Langlands correspondence. Then
\begin{equation*}
   \gamma_{\mathcal{A}}(s,\pi,\psi) = \gamma_{F}^{\rm Gal}(s,{}^\otimes{\rm I}(\sigma),\psi).
\end{equation*}
}

\bigskip

\noindent{\bf 3.4.} \emph{Proof of theorem.} Every irreducible representation $\pi$ of ${\rm GL}_n(F)$ can be written as the Langlands' quotient of a representation induced from quasi-tempered data (see \S~4.2 for more details). Notice that the Langlands quotient corresponds to taking the sum of the corresponding Weil-Deligne representations. Now, from properties (v) and (vi), we can assume that $\pi$ is cuspidal. We proceed by induction on $n \geq 1$, where the case $n=1$ is property~(iii).

{\sc Case 1}: $n>1$, $E/F$ tame, $\pi$ cuspidal and tame (thus $\pi$ of level zero). Let $\sigma$ be the transfer of $\pi$ under the local Langlands correspondence. Then $\sigma$ is irreducible and is given by
\begin{equation*}
   \sigma = {\rm Ind}_{\mathcal{W}_{E'}}^{\mathcal{W}_E}(\chi),
\end{equation*}
where $E'/E$ is an unramified extension of degree $n$ and $\chi : \mathcal{W}_{E'} \rightarrow \mathbb{C}^\times$ is a tame character. By class field theory, $\chi$ is identified with a character $\chi : E'^{\times} \rightarrow \mathbb{C}^\times$. Moreover, $\chi$ restricted to $U_{E'}$ is obtained from a regular character of $\mathds{k}_{E'}^\times$, i.e., a character with trivial stabilizer in ${\rm Gal}(\mathds{k}_{E'}/\mathds{k}_E)$. (Notation: given a local field $F$, its residue field is denoted by $\mathds{k}_F$; given a global function field $k$, its field of constants is denoted by $\mathds{k}_k$). 

Let $k = \mathds{k}_F(t)$ and let $v_0$ denote the place of $k$ given by the prime ideal $(t)$. Let $K/k$ be a separable quadratic extension such that there is unique place $w_0$ above $v_0$ and an isomorphism $\eta_0 : E/F \simeq K_{w_0}/k_{v_0}$. Such an extension has the property $\mathds{k}_K \simeq \mathds{k}_E$. Let $\mathds{k}_n$ be a degree $n$ extension of $\mathds{k}_K$, then the constant field extension $K' = \mathds{k}_n \cdot K$ is a cyclic extension of degree $n$, unramified everywhere. Then $K'$ has only one place $w_0'$ above $w_0$ and we can obtain an isomorphism $E'/E \simeq K_{w_0'}/K_{w_0}$ extending $\eta_0$. Let $w$ be a place of $K$ that splits completely in $K'/K$. (Notice that a place $w$ of $K$ splits completely in $K'/K$ if $n$ divides the degree of $w$). Let $S' = \left\{ w_0', w_1', \ldots, w_n' \right\}$, where $w_1', \ldots, w_n'$, are the places of $K'$ above $w$. We can now proceed as in \S~2.3 of \cite{hl} and construct a character
\begin{equation*}
   \xi : \prod U_{K_w'} \rightarrow \mathbb{C}^\times,
\end{equation*}
where the product ranges over all places $w'$ of $K'$, such that:
\begin{enumerate}
   \item $\xi_{w_0'}$ corresponds to $\chi \vert_{U_{E'}}$ via the isomorphism $K_{w_0'}' \simeq E'$;
   \item $\xi_{w'} = 1$ if $w' \notin S'$, and
   \item $\xi \vert_{\mathds{k}_{K'}^\times} = 1$.
\end{enumerate}
Then $\xi$ further extends to a gr\"ossencharacter
\begin{equation*}
   \tilde{\xi} : K'^\times \backslash \mathbb{A}_{K'}^\times \rightarrow \mathbb{C}^\times.
\end{equation*}
After globally twisting by a power of the norm in $\mathbb{A}_{K'}^\times$, we can assume that $\tilde{\xi}_{w_0'}$ corresponds to $\chi$ exactly. Notice that $\tilde{\xi}_{w'}$ will be unramified for $w' \notin S$.

A gr\"ossencharacter $\tilde{\xi}$ as above gives a character of $\mathcal{W}_{K'}$ via global class field theory, also denoted by $\tilde{\xi}$. Consider
\begin{equation*}
   R = {\rm Ind}_{\mathcal{W}_{K'}}^{\mathcal{W}_K} \tilde{\xi}.
\end{equation*}
Then $R_{w_0}$ corresponds to $\sigma$ via the isomorphism $\eta_0 : E \simeq K_{w_0}$. Also, $R_v$ will be reducible for all places $v$ of $K$ distinct from $w_0$; indeed, $R_v$ is unramified for $v \notin \left\{ w_0, w \right\}$, and $R_w$ is a sum of characters because $w$ is split in $K'/K$.

   Let $\ell$, $\ell \neq p$, be a fixed prime number, and let $\iota : \overline{\mathbb{Q}}_\ell \tilde{\, \rightarrow} \, \mathbb{C}$ be a fixed field isomorphism. Then, $R$ gives rise to a continuous degree $n$ $\overline{\mathbb{Q}}_\ell$-representation $\Sigma$ of $\mathcal{W}_K$. The global Langlands correspondence, proved in \cite{laff}, gives a cuspidal automorphic representation $\Pi = \Pi(\Sigma)$. By the local Langlands correspondence of \cite{lrs}, $\Pi_v$ corresponds to the Frobenius semisimplification of $\Sigma_v$.

%By construction: $\pi$ corresponds to $\Pi_{v_0}$, $\Pi_v$ is unramified for $v \notin \left\{ v_0, w \right\}$, and $\Pi_w$ is a principal series representation.

Choose a global character $\Psi = \otimes \Psi_v$ of $\mathbb{A}_k / k$ such that each $\Psi_v$ is non-trivial. Property~(iv) allows us to assume $\psi = \Psi_{w_0}$. The global functional equation then gives
\begin{equation} \label{eq:3.1}
   \prod_{v \in S} \gamma_{\mathcal{A}}(s,\Pi_v,\Psi_v) = \prod_{v \in S} \gamma_{k_v}^{\rm Gal}(s,{}^\otimes{\rm I}(\Sigma_v),\Psi_v),
\end{equation}
where $S$ is a finite set of places of $k$ containing $v_0$ and $u$, $w \vert u$. But for $v$ distinct from $v_0$ and $u$, $\Pi_v$ is unramified. Also, $\Pi_u$ is a (possibly ramified) principal series. Hence, properties~(iii) and (v) give
\begin{equation} \label{eq:3.2}
   \gamma_{K_v/k_v}(s,\Pi_v,r_{\mathcal{A}},\Psi_v) = \gamma_{k_v}^{\rm Gal}(s,{}^\otimes{\rm I}(\Sigma_v),\Psi_v), v \in S - \left\{ v_0 \right\}.
\end{equation}
Then (\ref{eq:3.1}) and (\ref{eq:3.2}) give equality at $v_0$. Notice that property~(i), used with the isomorphism $\eta_0: E/F \simeq K_{w_0}/k_{v_0}$, together with property~(ii) give the equality
\begin{equation*}
   \gamma_{E/F}(s,\pi,r_{\mathcal{A}},\psi) = \gamma_F^{\rm Gal}(s,{}^\otimes{\rm I}(\sigma),\psi),
\end{equation*}
in the case of $\psi = \Psi_{v_0} \circ \eta_0\vert_F$. And, property~(iv) gives the equality for any non-trivial $\psi$.

{\sc Case 2}: Let $n>1$, $E/F$ general and $\pi$ cuspidal, but not necessarily of level zero. Let $\sigma = \sigma(\pi)$ be the corresponding irreducible Weil-Deligne representation. Also, twisting $\pi$ by an unramified character, if necessary, gives a representation $\pi'$ with central character of finite order. The corresponding twist $\sigma'$ of $\sigma$ extends to a representation of ${\rm Gal}(\overline{F}/E)$, still denoted by $\sigma'$. Let $E'$ be the kernel of $\sigma'$, which is a finite galois extension of $E$, and let $\tilde{E}$ be the galois closure of $E'/F$ in $\overline{F}$. Consider, as above, the global function field $k = \mathds{k}_F(t)$ with its place $v_0$ and isomorphism $k_{v_0} \simeq F$. By the results of Gabber and Katz \cite{katz}, we can choose a finite galois extension $\tilde{K}$ of $k$ such that: it is tame at the place at infinity, denoted by $\infty$; it is unramified at all places $v \notin \left\{ v_0, \infty \right\}$; and, there is an isomorphism $\tilde{E}/F \simeq \tilde{K}_{\tilde{w}_0}/k_{v_0}$, where $\tilde{w}_0$ is the place of $\tilde{K}$ above $v_0$. Moreover, there is an isomorphism of ${\rm Gal}(\tilde{K}/k)$ onto ${\rm Gal}(\tilde{E}/F)$. We take for $K$ the quadratic separable extension of $k$ corresponding to $E$ via this isomorphism. Then, $\sigma'$ gives rise to a representation $\Sigma'$ of ${\rm Gal}(\tilde{K}/k)$. After twisting $\Sigma'$ by a power of the absolute value character, if necessary, we obtain a representation $\Sigma$ such that: $\Sigma_{w_0}$ corresponds to $\sigma$, where $w_0$ is the place of $K$ above $v_0$; it is tame at the place above $\infty$; and, it is unramified elsewhere.

By the global Langlands correspondence of \cite{laff}, there exists an irreducible cuspidal automorphic representation $\Pi$ corresponding to $\Sigma$. For every place $v$ of $K$, the local components $\Pi_v$ and $\Sigma_v$ correspond to each other via the local Langlands correspondence \cite{lrs}. Applying the same method as above, we get equality~(\ref{eq:3.1}) with $S = \left\{ v_0, \infty \right\}$. However, equality~(\ref{eq:3.2}) for $v = \infty$ is given by {\sc Case}~1, so that we obtain the desired equality at $v_0$. Finally, this is transferred to our given local situation as in the previous case.

\section{Local properties, $L$-functions and root numbers}

\noindent{\bf 4.1. Unitary groups and Asai $\gamma$-factors.} Let ${\bf G} = {\rm U}(2n)$ be a quasi-split unitary group with respect to $E/F$. For definiteness, we choose $\bf G$ to be defined via the hermitian form
\begin{equation*}
   h(x,y) = \sum_{i=1}^{n} \bar{x}_i y_{2n+1-i} - \sum_{i=1}^{n} \bar{x}_{2n+1-i}y_i
\end{equation*}
on the vector space $V = E^{2n}$. Let $J_n$ be the $n \times n$ matrix with $ij$-entries $\delta_{i,n-j+1}$. Then, take ${\bf M} \simeq {\rm Res}_{E/F}{\rm GL}_n$ to be the Siegel Levi subgroup of ${\bf G}$ whose group of rational points over $F$ is given by
\begin{equation*}
   M = \left\{ \left( \begin{array}{cc} g & 0 \\ 0 & J_n {}^t\bar{g}^{-1} J_n \end{array} \right) \vert g \in {\rm GL}_n(E) \right\}.
\end{equation*}

Let ${\bf B}$ be the borel subgroup of ${\rm GL}_n$ with unipotent radical ${\bf N}$ whose group of rational points $N$ over $E$ consists of upper triangular matrices. Given a non-trivial character $\psi : F \rightarrow \mathbb{C}$, it is extended to a character of $N$, also denoted by $\psi$, by setting $\psi(x) = \psi(x_{1,2} + \cdots + x_{n-1,n})$, $(x) = (x_{i,j}) \in N$.

Given $(E/F,\pi,\psi) \in \mathscr{A}_{\rm quad}(p)$ generic, the definition of $\gamma$-factors is given via the local coefficient $C_\chi(s,\pi,w_0)$ in \S~6.4 of \cite{luis}. We can assume that $\pi$ is $\psi$-generic. Let $\bf T$ be the maximal torus of $\bf G$. For $a \in F^\times$, let 
\begin{equation*}
   t={\rm diag}(a^{-(n-1)}\alpha, a^{-(n-2)}\alpha, \ldots, a\alpha, \alpha, \beta, a^{-1}\beta, \ldots, a^{n-2}\beta, a^{n-1}\beta) \in {\bf T}(\overline{F}),
\end{equation*}
where $\alpha \beta^{-1} = \bar{\alpha}^{-1}\bar{\beta} = a \in F^\times$. Then $w_0(t^{-1})t$ is an element of the center of $M$. Let $\pi_t$ be given by $\pi_t(x) = \pi(t^{-1}xt)$. The character $\psi_t$ of $N$ given by $\psi_t(x) = \psi(t^{-1}xt)$, is obtained from the non-trivial character $\psi^a$ of $F$ as above. The representation $\pi_t$ is then $\psi_t$ generic.  Using equation~(6.1) of \S~6.2 of \cite{luis} together with the definition we explicitly get
\begin{align*}
   \gamma_{E/F}(s,\pi,r_\mathcal{A},\psi^a) 	& = C_{\psi^a}'(s,\pi_t,w_0) \\
   									& = \omega_\pi(a)^n \left| a \right|_F^{n^2(s-\frac{1}{2})} C_{\psi}(s,\pi,w_0) \\
									& = \omega_\pi(a)^n \left| a \right|_F^{n^2(s-\frac{1}{2})} \gamma_{E/F}(s,\pi,r_\mathcal{A},\psi).
\end{align*}

\noindent{\bf 4.2. Local Asai factors.} Let us recall the definition of Asai local factors for generic $\pi$, and then extend it to any smooth irreducible represenation $\pi$ of ${\rm GL}_n(E)$. Since the local factors we are studying arise from representations of ${\rm GL}_n$, we extend this definition to $(E/F,\pi,\psi) \in \mathscr{A}_{\rm quad}(p)$ where $\pi$ is any smooth representation of ${\rm GL}_n(E)$.

Let $(E/F,\pi,\psi) \in \mathscr{A}_{\rm quad}(p)$. Let us first assume that $\pi$ is tempered, then $\pi$ is generic \cite{ze}. Let $P_{\pi}(t)$ be the unique polynomial satisfying $P_\pi(0) = 1$ and such that $P_\pi(q^{-s})$ is the numerator of $\gamma_{E/F}(s,\pi,r_{\mathcal{A}},\psi)$. Then
\begin{equation*}
   L(s,\pi,r_{\mathcal{A}}) = P_\pi(q^{-s})^{-1}.
\end{equation*}
Because $\pi$ is tempered, $L(s,\pi,r_{\mathcal{A}})$ is holomorphic for $\Re(s) > 0$. If $\pi$ is obtained by unitary parabolic induction from $\pi_1 \otimes \cdots \otimes \pi_d$, where each $\pi_i$ is tempered, then multiplicativity of $\gamma$-factors gives multiplicativity of $L$-functions:
\begin{equation*}
   L(s,\pi,r_{\mathcal{A}}) = \prod_{i=1}^d L(s,\pi_i,r_{\mathcal{A}}) \prod_{i<j} L(s,\pi_i \times \pi_j^{\rm conj}).
\end{equation*}
The local $\varepsilon$-factor is defined to satisfy the relation:
\begin{equation*} \label{rootnumber}
   \gamma_{E/F}(s,\pi,r_{\mathcal{A}},\psi) = \varepsilon_{E/F}(s,\pi,r_{\mathcal{A}},\psi) \dfrac{L(1-s,\widetilde{\pi},r_{\mathcal{A}})}{L(s,\pi,r_{\mathcal{A}})}.
\end{equation*}

Given $(E/F,\pi,\psi) \in \mathscr{A}_{\rm quad}(p)$ in general, we can use Langlands classification to write $\pi$ as the Langlands quotient of a representation
\begin{equation*}
   \tau = {\rm ind}_{P'}^{{\rm GL}_n(E)} (\tau_1 \otimes \cdots \otimes \tau_k).
\end{equation*}
Here, ${\bf P}'$ is the parabolic subgroup of ${\rm GL}_n(E)$ with Levi subgroup $\prod_{i=1}^d {\rm GL}_{n_i}(E)$; each $\tau_i$ is a quasi-tempered representation of ${\rm GL}_{n_i}(E)$ of the form $\left| {\rm det} (\cdot) \right|_E^{t_i} \tau_{i,0}$; and, $t_1 > \cdots > t_d$ is a sequence of real numbers. Each $\tau_{i,0}$ is tempered, and we define Asai $\gamma$-factors in general by setting
\begin{equation*}
   \gamma_{E/F}(s,\pi,r_{\mathcal{A}},\psi) = \prod_{i=1}^d \gamma_{E/F}(s + 2t_i, \tau_{i,0},r_{\mathcal{A}},\psi) 
   								     \prod_{i<j} \gamma(s+t_i+t_j,\tau_{i,0} \times \tau_{j,0}^{\rm conj},\psi_E),
\end{equation*}
where $\psi_E = \psi \circ {\rm Tr}_{E/F}$. For generic triples the definition is agrees with that of \cite{luis}, where the multiplicativity property of Asai $\gamma$-factors takes the following form:
\begin{itemize}
   \item[(v)'] (Multiplicativity). \emph{For $i = 1, \ldots, d$, let $(E/F,\pi_i,\psi) \in \mathscr{A}_{\rm quad}(p)$ be generic of degree $n_i$. Let $n = n_1 + \cdots + n_d$ and let $\pi$ be the unique generic constituent of the representation of ${\rm GL}_n(E)$, obtained via unitary parabolic induction from $\pi_1 \otimes \cdots \otimes \pi_d$. Then
   \begin{equation*}
      \gamma_{E/F}(s,\pi,r_{\mathcal{A}},\psi) = \prod_{i=1}^d \gamma_{E/F}(s,\pi_i,r_{\mathcal{A},\psi}) 
      									\prod_{i<j} \gamma_E(s,\pi_i \times \pi_j^{\rm conj},\psi_E).
   \end{equation*}
   }
\end{itemize}
It is now easy to check that properties (i) through (vii) are satisfied by Asai $\gamma$-factors. Any rule for generic tripes which satisfies these properties is unique and is given by the corresponding Galois factors. Actually, the proof of Theorem~3.3 does not require property~(vi) in this case. More explicitly, let $(E/F,\pi,\psi) \in \mathcal{A}_{\rm quad}(p)$ be generic, then
\begin{equation*}
   \gamma_{\mathcal{A}}(s,\pi,\psi) = \gamma_F^{\rm Gal}(s,{}^{\otimes}{\rm I}(\sigma),\psi).
\end{equation*}
In addition to the properties in the characterization, the local functional equation of $\gamma_F^{\rm Gal}$ can be immediately translated into a property of $\gamma_{\mathcal{A}}$ via the main theorem:
\begin{itemize}
   \item[(viii)] (Local functional equation). \emph{Let $(E/F,\pi,\psi) \in \mathscr{A}_{\rm quad}(p)$, then
      \begin{equation*}
         \gamma_{\mathcal{A}}(s,\pi,\psi) \gamma_{\mathcal{A}}(1-s,\widetilde{\pi},\overline{\psi}) = 1.
      \end{equation*}
      }
\end{itemize}

Now, let $(E/F,\pi,\psi) \in \mathscr{A}(p)$ and write $\pi$ as a Langlands quotient $\pi = {\rm J}(\tau)$, with $\tau$ as above. Asai $L$-functions and $\varepsilon$-factors are defined accordingly:
\begin{align*}
   L(s,\pi,r_{\mathcal{A}}) & = \prod_{i=1}^{d} L(s + 2t_i,\tau_{i,0},r_{\mathcal{A}}) 
   					  \prod_{i<j} L(s + t_i + t_j,\tau_{i,0} \times \tau_{j,0}^{\rm conj}), \\
   \varepsilon_{E/F}(s,\pi,r_{\mathcal{A}},\psi) & = \prod_{i=1}^d \varepsilon_{E/F}(s + 2t_i, \tau_{i,0},r_{\mathcal{A}},\psi)
   									    \prod_{i<j} \varepsilon(s+t_i+t_j,\tau_{i,0} \times \tau_{j,0}^{\rm conj},\psi_E).
\end{align*}

\bigskip

\noindent{\bf 4.3. Theorem.} \emph{If $(E/F,\pi,\psi) \in \mathscr{A}_{\rm quad}(p)$, let $(E/F,\sigma,\psi) \in \mathscr{G}_{\rm quad}(p)$ be the triple obtained by taking $\sigma = \sigma(\pi)$ to be the transfer of $\pi$ under the local Langlands correspondence. Then
\begin{align*}
   L(s,\pi,r_\mathcal{A}) & = L(s,{}^\otimes{\rm I}(\sigma)) \\
   \varepsilon_{E/F}(s,\pi,r_\mathcal{A},\psi) & = \varepsilon(s,{}^\otimes{\rm I}(\sigma),\psi).
\end{align*}
}

\noindent{} \emph{Proof.} The equality
\begin{equation*}
   \gamma_{E/F}(s,\pi,r_\mathcal{A},\psi) = \gamma_F^{\rm Gal}(s,{}^\otimes{\rm I}(\sigma),\psi)
\end{equation*}
follows from Theorem~3.3 for tempered $\pi$, since tempered representations are generic. The local Langlands correspondence sends $\pi$ tempered to a sum of indecomposable Weil-Deligne representations with unitary central character. Then ${}^\otimes{\rm I}(\sigma)$ is also a such a direct sum and consequently the $L$-function $L(s,{}^\otimes{\rm I}(\sigma))$ is holomorphic for $\Re(s)>0$. Hence, there is no cancelation in equation~(4.1) for $0 < \Re(s) < 1$ and $L(s,{}^\otimes{\rm I}(\sigma))$ is exactly the denominator of $\gamma_{E/F}(s,\pi,r_\mathcal{A},\psi)$. For $\pi$ in general, the Langlands quotient corresponds to a direct sum of indecomposable representations. The equalities of the theorem now follow, since local $L$-functions and $\varepsilon$-factors for $\mathscr{G}_{\rm quad}(p)$ are compatible with the above definition (see \S~2.9 of \cite{he02}).

\bigskip

\noindent{\bf 4.4. Stability of local factors.} Let us briefly recall the stability property of local factors for Weil-Deligne representations \cite{d,dh}. Let $\sigma$ be a Weil-Deligne representation. Let $\eta$ be a character of $F^\times$ of level $k$, for $k$ sufficiently large (depending on $\sigma$). Take an element $c = c(\eta,\psi) \in F^\times$ such that $\psi(cx) = \eta(1+x)$ for $x \in \mathfrak{p}^{\left[ k/2 \right] + 1}$. Then
\begin{equation*}
   \varepsilon(s,\sigma \otimes \eta,\psi) = {\rm det}(\sigma(c))^{-1} \varepsilon(s,\eta,\psi)^{{\rm dim}\sigma},
\end{equation*}
\begin{equation*}
   L(s,\sigma \otimes \eta) = L(s,\widetilde{\sigma} \otimes \eta^{-1}) = 1.
\end{equation*}

Because of this, the next property is now a corollary to Theorem~4.4.

\bigskip

\noindent{\bf Corollary} (Stability). \emph{Let $(E/F,\pi_i,\psi) \in \mathscr{A}_{\rm quad}(p)$, $i = 1,2$, be of the same degree. Assume that $\pi_1$ and $\pi_2$ have the same central character. If $\eta$ is a character of $F^\times$ such that $\eta \cdot \eta^{\rm conj}$ is sufficiently highly ramified, then
\begin{equation*}
   \varepsilon_{E/F}(s,\eta \cdot \pi_1,\psi) = \varepsilon_{E/F}(s,\eta \cdot \pi_2,\psi)
\end{equation*}
\begin{equation*}
   L(s,\eta \cdot \pi_1) = L(s,\eta^{-1} \cdot \widetilde{\pi}_2) = 1.
\end{equation*}
In terms of $\gamma$-factors
\begin{equation*}
   \gamma_{\mathcal{A}}(s,\eta \cdot \pi_1,\psi) = \gamma_{\mathcal{A}}(s,\eta \cdot \pi_2,\psi).
\end{equation*}
}

\bigskip

\noindent{\bf 4.5. A note on exterior and symmetric square $\gamma$-factors.} A similar stability property also holds for local factors corresponding to exterior and symmetric square $L$-functions. To be more precise, we use the notation of \cite{hl}: Let $(F,\psi,\pi_i) \in \mathscr{L}(p)$, $i = 1,2$, be of degree $n$. Assume that $\pi_1$ and $\pi_2$ have the same central character. If $\eta$ is a character of $F^\times$ such that $\eta^2$ is sufficiently highly ramified, then
\begin{equation*}
   \gamma_F(s,\eta \cdot \pi_1, r_n,\psi) = \gamma_{F}(s,\eta \cdot \pi_2,r_n,\psi),
\end{equation*}
where $r_n$ is either $\wedge^2 \rho_n$ or ${\rm Sym}^2 \rho_n$. Furthermore, notice that the discussion of \S~4.2 can be adapted to this situation. After incorporating twists by unramified characters and Langlands classification to the characterization of [loc. cit.], this result holds in the general setting of smooth irreducible representations $\pi_i$ of ${\rm GL}_n(F)$.

\end{document}